\documentclass[journal,10pt,twocolumn]{IEEEtran}
\usepackage[fleqn]{amsmath}
\usepackage{graphicx,subfigure,amsthm,array,amsfonts,amssymb,bm,mathrsfs,extarrows}
\usepackage{cite}
\usepackage{color}
\usepackage{xcolor}
\theoremstyle{remark}

\newtheorem{remark}{Remark}
\newtheorem{theorem}{Theorem}

\newtheorem{lemma}{Lemma}
\newtheorem{proposition}{Proposition}

\usepackage[switch]{lineno}
\allowdisplaybreaks

\begin{document}
\title{Formation Control with Unknown Directions and General Coupling Coefficients}

\author{Zhen Li,
        Yang Tang \emph{Member, IEEE},
        Yongqing Fan,
        and Tingwen Huang, \emph{Fellow, IEEE}
\thanks{This work was supported in part by the National Natural Science Foundation of China under Grant 62003262, the Programme of Introducing Talents of Discipline to Universities (the 111 Project) under Grant B17017, and the National Priority Research Project from the Qatar National Research Fund under Grant NPRP 9-166-1-031. \emph{(Corresponding
authors: Yang Tang and Tingwen Huang.)}}
\thanks{Zhen Li is with the School of Automation, Xi-an University of Posts $\&$ Telecommunications, Xi-an 710121, China (e-mail: lzqz12@hotmail.com). }
\thanks{Yang Tang is with the Key Laboratory of Smart Manufacturing in Energy Chemical Process, East China University of Science and Technology, Shanghai 200237, China (e-mail: tangtany@gmail.com). }
\thanks{Yongqing Fan is the School of Automation, Xi-an University of Posts $\&$ Telecommunications, Xi-an 710121, China (e-mail: fanyongqing@xupt.edu.cn).}
\thanks{Tingwen Huang is with the Science Program, Texas A$\&$M University, Doha 23874, Qatar (e-mail: tingwen.huang@qatar.tamu.edu).}
 }
\maketitle

\begin{abstract}
Generally, the normal displacement-based formation control has a sensing mode that requires the agent not only to have certain knowledge of its direction, but also to gather its local information characterized by nonnegative coupling coefficients. However, the direction may be unknown in the sensing processes, and the coupling coefficients may also involve negative ones due to some circumstances. This paper introduces these phenomena into a class of displacement-based formation control problem. Then, a geometric approach have been employed to overcome the difficulty of analysis on the introduced phenomena. The purpose of this approach is to construct some convex polytopes for containing the effects caused by the unknown direction, and to analyze the non-convexity by admitting the negative coupling coefficients in a certain range. Under the actions of these phenomena, the constructed polytopes are shown to be invariant in view of the contractive set method. It means that the convergence of formation shape can be guaranteed. Subsequently, an example is given to examine the applicability of derived result.
\end{abstract}

\begin{IEEEkeywords}
Unknown direction, general coupling coefficients, formation control, polytope
\end{IEEEkeywords}
\IEEEpeerreviewmaketitle

\section{Introduction}\label{S1}
Recently, formation control problem has been getting more and more research attentions in many branches of engineering, which has also offered theoretical support for practical applications, e.g., flying formation of unmanned aerial vehicles \cite{Bayezit}, target enclosing \cite{Zheng}, and source seeking \cite{Han2014}. The main purpose of the considered problem is to construct a mechanism, such that the agents could automatically tend to a specific geometric shape. For the sake of realizing this purpose, each agent needs to regulate the current state according to the non-centralized information from its neighbors. This process can be described by different control mechanisms that are usually classified into three categories: position-based, displacement-based, and distance-based methods \cite{Oh2015}. Among these control mechanisms, the displacement-based method has been intensively considered into multi-agent system due to high efficiency and simple structure \cite{Verginis,Aranda,Montijano,Thunberg,Marina,Dong2019,Zhang2021,Cao,Wu2013,Tang2020}.

For the displacement-based method, one of the problems to be considered is the sensing ability of the agent \cite{Oh2015,Tang2013}. This ability is often expressed by the sensing mode dependent on the actual sensor, and further determines the type of controlled variables. Specifically, in order to realize the desired formation shape, the agent senses and drives the relative displacements based on its neighbors as an active way. There is a precondition that, for each agent, the relative displacements are assumed to be measured according to its local coordinate frame, whereas the global coordinate frame is not necessary to be completely known for all the agents. However, the direction between the local and the global coordinate frames has to be aligned, otherwise the actual formation shape may be distorted by the misalignment as time evolves \cite{Meng2016,Li2019,Tang2022}. For the latter, there have been some results that consider the compensation control by utilizing the position filtering for reconstructing the formation control strategy \cite{Meng2016}. Meanwhile, some other works have studied similar problems in view of the relative direction measurement \cite{Tran,Lee,Li2019,Jin}.

Actually, for each agent, the misalignment between the local and the global coordinate frames results in a rotation angle. Although some existing works have shown that this angle may be irrelevant to the system performance, it is mainly related to the distance-based method \cite{Lin,Zhao}. For the distance-based method, the convergence of the desired formation shape has been replaced by the potential of system structures and configuration, such that the rotation may not affect the convergence under the graph rigidity theory \cite{Han2016,Li2018}. By contrast, it usually requires to design an additional control strategy for aligning the direction under the displacement-based method \cite{Wang2021}. This further means that the agent has to collect the enough direction information for the alignment. However, the measurement of direction is not always reliable. There are some examples in the actual actuators and sensors, e.g., errors accumulation of the inertial navigation, and local anomalies of the geomagnetic fields \cite{Yu2016,Meng2016}. Therefore, when the rotation angle is unreliable or even unknown, it is significant to consider a control strategy in order to achieve the desired formation shape.

Another problem in the displacement-based method involves the interaction topology \cite{Oh2015}. The interaction topology is attributed to the sensing ability of the agents, whose edges are often assumed to be characterized by nonnegative coupling coefficients \cite{Dong2019,Zhang2021}. These nonnegative coupling coefficients usually express the cooperative interactions among the agents. However, general coupling coefficients may emerge from some cases, e.g., some auxiliary systems \cite{Liu2020}, and the social networks \cite{Altafini}. When the relationship between the interacting agents is characterized by the general coupling coefficients, both the cooperative and competitive interactions may arise, where the negative coupling coefficients represent the competitive interactions. The competitive interactions in some related works have been utilized to describe a disagreement tendency that is known as the formation and separation problems \cite{Li2022}. It needs to be emphasized that, in these works \cite{Li2022}, the desired formation shape evolves into two parts separated from each other, which is eventually corrupted by the competitive interactions. Therefore, how to promote and keep the desired formation shape under the coexistence of the cooperative and competitive interactions is still an open problem.

On the other hand, the usual method for investigating formation control problems is to transform the system into a convergence problem. Generally, the Lyapunov function plays an important factor for analyzing this convergence problem, where the convergence rate can be solved by some inequalities satisfying the existence of the Lyapunov function \cite{Marina,Dong2019,Zhang2021,Liu2022}. However, as mentioned in \cite{Olshevsky,Chen2016,Matei},
the existence of the common Lyapunov function may not be guaranteed for some time-varying cases, not to mention solving the convergence rate. Recently, by utilizing the convex polytope, an alternative method has been proposed for dealing with the cooperative problems \cite{Chen2013}. The main line of this method is to treat the formation control protocol as a convex combination, whose convergence can be solved by the principle of contractive set. It should be noted that the misalignment coordinate frames and the general coupling coefficients may cause the non-convexity, since the special orthogonal group $\mathrm{SO}(n)$ and the Euclidean space are non-homeomorphic spaces\cite{Saunderson}, and the general coupling coefficients have negative terms. These imply that the convex polytope may hard to be directly applied. In this case, it is necessary to extend an alternative way for investigating the formation control problem under the misalignment coordinate frames and the general coupling coefficients.

Motivated by the aforementioned observation, a kind of the displacement-based formation control problem has been considered in this paper, and the corresponding analysis method are further developed. When designing the control protocol, neither each agent is required to have the ability for sensing the direction, nor the coupling coefficients between any pair of the agents are nonnegative. For the sake of understanding the non-convexity caused by these cases, several convex polytopes are introduced, and their properties are then discussed. The obtained results show that, although the convexity of these polytopes can be ensured under the actions of the unknown directions and the general coupling coefficients, the contractility may be collapsed. Then, the polar of the polytopes is considered for describing a larger convex region, such that the polytopes are compatible with the behaviors caused by the unknown directions and the general coupling coefficients. The contracting mapping principle is applied to designing some parameters in the case that this convex region can be contractive. Based on the obtained results, the convergence rate is solved, which means that the desired formation shape can be finally realized.

The contributions of obtained results are characterized by the following points: 1) Compared with the normal displacement-based formation control problems \cite{Dong2019,Zhang2021}, the alignment between the local and the global coordinate frames is not require for each agent. Different from the formation control problems \cite{Tran,Lee,Meng2016,Li2019}, the direction of the agent in this paper does not need to be estimated or measured during the control process, which can even be completely unknown. 2) When designing the control protocol, the general coupling coefficients are introduced. In this case, the coupling coefficients are not only nonnegative, but also allowed to be negative within a certain range. Different from the negative coupling coefficients in the formation and separation problems \cite{Yan2022,Li2022}, the desired formation shape can still be guaranteed under the action of the general coupling coefficients. 3) Inspired by \cite{Chen2013}, the convex polytope method is developed to study the displacement-based formation control for considering the unknown directions and the general coupling coefficients. In this case, some additional restrictions of the Lyapunov function are removed, e.g., the static interactions or the existence of system eigenvalues \cite{Zhang2021,Cao}. Meanwhile, the method in this paper also extends some linear relevant results \cite{Chen2013,Chen2016,Liu2020} onto $\mathrm{SO}(2)$.

A outline of the remainder of this paper is as follows. In Section II the signed graph, the convex polytope and some preliminaries are expressed, then the concerned formation control problem is introduced. In Section III, main results are provided by utilizing some properties of the convex polytope. Section IV gives a simulation example, and then the conclusions are shown in Section V. Section VI includes the technical proofs of some results in Section III.

\emph{Notation}: $\mathbb{R}^{N\times M}$ and $\mathbb{R}^{n}$ are utilized to express the set of $N\times M$ real matrices and $n$-dimensional real vector space, respectively. $\textbf{0}_{N}$ and $\textbf{1}_{N}$ stand for a unit vector and a null vector in $ \mathbb{R}^{N}$, respectively. $e^{i}_{N}$ is the basis vector in $\mathbb{R}^{N}$, whose $i$-th entry is $1$ and all rests are $0$. $O_{N}$ and $I_{N}$ are a null matrix and an identity matrix in $ \mathbb{R}^{N\times N}$, respectively. $\mathbb{N}$ and $\mathbb{N}^{+}$ stand for the set of natural and positive natural numbers, respectively. $\langle x,y\rangle$ is the inner product for any $x,y\in\mathbb{R}^{N}$. $\|x\|_{p}$ represents the $p$-matrix norm for a matrix $x$, or the $\ell_{p}$ vector norm for a vector $x$. $\otimes$ is the Kronecker product. $\mathrm{SO}(2)=\{R\in\mathbb{R}^{2\times 2}|R^{T}R=RR^{T}=I_{2},\det(R)=1\}$ is the special orthogonal group defined on $\mathbb{R}^{2\times 2}$, where $\det (R)$ is the determinant of $R$.

\section{preliminaries and problem formulation}\label{S2}

\subsection{Some preliminaries}\label{S21}

Denote a pair $(\mathcal{V}_{\mathcal{G}}, \mathcal{E})$ as a directed graph $\mathcal{G}$, which consists of the vertex set $\mathcal{V}_{\mathcal{G}}=\{1,\ldots,N\}$ and the edge set $\mathcal{E}\subseteq\mathcal{V}_{\mathcal{G}}\times\mathcal{V}_{\mathcal{G}}$. For a signed graph $\mathcal{G}$, the edge set $\mathcal{E}=\mathcal{E}^{+}\cup\mathcal{E}^{-}$ is generated by the positive edge set $\mathcal{E}^{+}$ and the negative edge set $\mathcal{E}^{-}$. The signed directed graph $\mathcal{G}$ is described by a \emph{weighted adjacency matrix} $\mathcal{A}=[a^{ij}]\in\mathbb{R}^{N\times N}$ whose entries satisfy: $a^{ij}\neq 0$ if $(j,i)\in\mathcal{E}$ holds for $i\neq j$, otherwise, $a^{ij}=0$, and particularly, $a^{ii}=0$. In view of the signed graph theory, for any distinct $i,j\in\mathcal{V}_{\mathcal{G}}$, the interaction between the vertexes $i$ and $j$ is \emph{cooperative} if $a^{ij}>0$ holds for $(j,i)\in\mathcal{E}^{+}$, and the interaction between the vertexes $i$ and $j$ is \emph{competitive} if $a^{ij}<0$ holds for $(j,i)\in\mathcal{E}^{-}$. For the graph $\mathcal{G}$, its Laplacian matrix $\mathcal{L}=[\ell^{ij}]\in\mathbb{R}^{N\times N}$ corresponding to the weighted adjacency matrix $\mathcal{A}$ satisfies: $\ell^{ij}=-a^{ij}$ for any $i\neq j$, otherwise, $\ell^{ii}=-\sum_{j\in\mathcal{N}^{i}}a^{ij}$, where the neighbor set of vertex $i$ is $\mathcal{N}^{i}=\{j\in\mathcal{V}_{\mathcal{G}} |  (i,j)\in\mathcal{E}\}$. A directed graph $\mathcal{G}$ is called \emph{neighbor shared} \cite{Chen2013}, if for any distinct $p,q \in\mathcal{V}_{\mathcal{G}}$, there exists $i\in\mathcal{V}_{\mathcal{G}}$ such that $(i,p)\in\mathcal{E}^{+}$ and $(i,q)\in\mathcal{E}^{+}$ hold for $p\neq i$ and $q\neq i$. Apparently, a graph must be connected, if it is neighbor shared.

A matrix $S=[s^{ij}]\in\mathbb{R}^{N\times N}$ is called \emph{nonnegative}, if $s^{ij}\geq0$ holds for all $i,j$. The nonnegative matrix $S$ is called \emph{stochastic}, if $\sum_{j=1}^{N}s^{ij}=1$ holds for each $i$. The matrix $S$ is called \emph{general stochastic}, if $\sum_{j=1}^{N}s^{ij}=1$ holds for each $i$, whose entries $s^{ij}$ are not necessarily nonnegative. For the general stochastic matrix $S$, its \emph{ergodic coefficient} $\eta(S)$ is defined by $\eta(S)=\max_{p,q}\sum_{i=1}^{N}|s^{pi}-s^{qi}|/2$. Denote a matrix set $\mathbb{S}_{\beta}$  \cite{Chen2016,Liu2020} as:
\begin{enumerate}
  \item $\sum_{j=1}^{N}s^{ij}=1$
  \item $\sum_{j=1}^{N}\min\{s^{ij},0\}\geq-\beta$ for each $i$,
\end{enumerate}
where $s^{ij}$ is the $ij$-entry of the matrix $S\in\mathbb{S}_{\beta}$, and $\beta\geq0$. $\mathbb{S}_{\beta}$ is obviously the set of general stochastic matrices if $\beta\neq 0$, and specifically, $\mathbb{S}_{0}$ is the set of stochastic matrices.

For a subset $\mathcal{V}$ of $ \mathbb{R}^{N}$, its \emph{convex hull} is the smallest convex set containing $\mathcal{V}$, i.e.,
\begin{align*}
\text{conv}\mathcal{V}=&\Big\{\sum_{i=1}^{\#\mathcal{V}}\lambda^{i}x^{i}\Big| \sum_{i=1}^{\#\mathcal{V}}\lambda^{i}=1,\lambda^{i}\in[0,1],x^{i}\in\mathcal{V}\Big\},
\end{align*}
where $\#\mathcal{V}$ is the associated cardinality. Specially, if $\#\mathcal{V}<\infty$, $\text{conv}\mathcal{V}$ expresses a $\mathcal{V}$-\emph{polytope} \cite{Ziegler}, denoted by $\mathcal{P}_{\mathcal{V}}$. The face of dimension $0$ is a \emph{vertex} of $\mathcal{P}_{\mathcal{V}}$. A point is called an \emph{extreme point} of $\mathcal{P}_{\mathcal{V}}$, if it can not express this point as a convex combination of at least two distinct points in $\mathcal{P}_{\mathcal{V}}$. From \cite{Lay}, a vertex of any compact convex set is automatically an extreme point, which means that the sets of all vertices and extreme points of $\mathcal{P}_{\mathcal{V}}$ are equivalent, defined as $\text{ext}\mathcal{P}_{\mathcal{V}}$. Moreover, according to the Minkowski-Weyl theorem \cite{Rockafellar}, $\mathcal{P}_{\mathcal{V}}$ has a nonempty intersection form of a sequence of closed halfspaces, denoted by an $\mathcal{H}$-\emph{polytope} $\mathcal{P}_{N}=\{x\in\mathbb{R}^{N}|Ax\leq z\}$ for some $A\in\mathbb{R}^{M\times N}$ and $z\in\mathbb{R}^{M}$ \cite{Ziegler}. For a matrix $\mathcal{M}\in\mathbb{R}^{M\times N}$, $\mathcal{M}\mathcal{P}=\{\mathcal{M}x|x\in\mathcal{P},\mathcal{P}\subseteq \mathbb{R}^{N}\}$ means a linear mapping of the polytope $\mathcal{P}$. For a polytope $\mathcal{P}\subseteq\mathbb{R}^{N}$, its polar is denoted by
\begin{align*}
\mathcal{P}^{\circ}=&\{x\in\mathbb{R}^{N}| \langle x,y\rangle\leq 1,\text{for all }y\in\mathcal{P}  \}.
\end{align*}

\subsection{Problem formulation}\label{S22}
For any $i\in\mathcal{V}_{\mathcal{G}}$, denote $\mathfrak{F}^{w}$ and $\mathfrak{F}^{i}$ as the subspaces of $\mathbb{R}^{2}$, where, for all agents, $\mathfrak{F}^{w}$ means the global coordinate frame (fixed), and for $i$-th agent, $\mathfrak{F}^{i}$ is the local coordinate frame. Consider a group of $N$ agents composed of first-order dynamics in discrete-time domain with respect to $\mathfrak{F}^{w}$, where the $i$-th agent is modeled as follows
\begin{align}
\Delta_{h_{k}} p^{i}_{k}= u^{i}_{k},\label{eq1}
\end{align}
where $\Delta_{h_{k}}$ is the forward-difference operator, $h_{k}=t_{k+1}-t_{k}$ ($k\in\mathbb{N}$) is a sampling interval, $t_{k}$ is the time instant, $p^{i}_{k}=[p_{k}^{i1},p_{k}^{i2}]^{T}\in \mathfrak{F}^{w}$ expresses the position vector at $t_{k}$, the position vector $p^{i}_{k}$ is initialized at $t_{0}$, and $u^{i}_{k}$ stands for the control protocol at $t_{k}$.

Generally, each agent drives the current relative positions (displacements) measured from the neighbors to approach the desired positions, namely, the desired formation shape. In order to realize this aim, all the agents should have some external measuring abilities with respect to $\mathfrak{F}^{w}$. However, due to the limitations on the types and capacities of sensors, the agent $i$ can only measure the relative positions with respect to $\mathfrak{F}^{i}$. Moreover, the desired positions mainly have two views. The first one is to regard the desired positions as some reference inputs with respect to $\mathfrak{F}^{w}$ \cite{Oh2015}. The second one treats the desired positions as analogous concepts of fixed landmarks in the navigation problems \cite{Cao}, which requires to sense the differences between the interacting agents and corresponding desired positions. In this paper, the desired positions are assumed to follow the second view, and are fixed with respect to $\mathfrak{F}^{w}$. When the agent $i$ senses the information with respect to $\mathfrak{F}^{i}$, the desired positions become time-varing due to the motion of the agent $i$. In this case, for the agent $i$ and its neighbor $j$, the sensed information is expressed by
\begin{align}
\xi^{ji}_{k}=p^{j,i}_{k}-p^{i,i}_{k}-(d^{j,i}_{k}-d^{i,i}_{k}).\label{eq2}
\end{align}
where $p^{j,i}_{k}\in\mathfrak{F}^{i}$ is the position vector of the agent $j$, $p^{i,i}_{k} $ is the origin of $\mathfrak{F}^{i}$, $d^{j,i}_{k}$ and $d^{i,i}_{k}\in\mathfrak{F}^{i}$ are the desired position vectors of the agents $j$ and $i$ at $t_{k}$, respectively.

However, the direction sensors may have an inaccurate result, such that the dynamics of each agent may not match $\mathfrak{F}^{w}$. It further implies that $\mathfrak{F}^{w}$ and $\mathfrak{F}^{i}$ may be misaligned. Therefore, this paper is assumed that each agent has a misaligned direction between $\mathfrak{F}^{w}$ and $\mathfrak{F}^{i}$, where the rotation angle induced by this misalignment is unknown. In order to describe this misalignment, a rotation matrix $R^{i}_{k}\in SO(2)$ is presented by
\begin{align}
R^{i}_{k}=\cos\theta_{k}^{i}I_{2}+\sin\theta_{k}^{i}J,\label{eq3}
\end{align}
where the matrix $J=\big[\begin{smallmatrix}
0&  -1\\
     1& 0
\end{smallmatrix}\big]$, and $\theta_{k}^{i}\in(-\pi/2,\pi/2]$ is a unknown rotation angle indicated by the misaligned direction between $\mathfrak{F}^{w}$ and $\mathfrak{F}^{i}$ at $t_k$.

Note that $\xi^{ji}_{k}$ can be considered as a special kind of output measurement with respect to $\mathfrak{F}^{i}$ rather than $\mathfrak{F}^{w}$. Before designing the control protocol $u^{i}_{k}$ in (\ref{eq1}), one needs to transform $\xi^{ji}_{k}$ into a measurement with respect to $\mathfrak{F}^{w}$ for the sake of facilitating the analysis. From (\ref{eq2}) and (\ref{eq3}), this measurement is expressed by
\begin{align}
\xi^{ji}_{k}=(R^{i}_{k})^{T}\big(p^{j}_{k}-p^{i}_{k}-(d^{j}-d^{i})\big),\label{eq4}
\end{align}
where $p^{j}_{k}-p^{i}_{k}=R^{i}_{k}(p^{j,i}_{k}-p^{i,i}_{k})$, $d^{j}-d^{i}=R^{i}_{k}(d^{j,i}_{k}-d^{i,i}_{k})$, and $d^{i}\in\mathfrak{F}^{w}$ means the fixed desired position vector of the agent $i$ at $t_{k}$.

Similar with the sensing mode of the agents $i$ and $j$ from (\ref{eq4}), gathering all the measurements in $\mathcal{N}^{i}_{k}$, the control protocol $u^{i}_{k}$ in (\ref{eq1}) is given as follows
\begin{align}
u^{i}_{k}=  &\sum_{j\in\mathcal{N}^{i}_{k}} a^{ij}_{k}(R^{i}_{k})^{T}\big(p^{j}_{k}-p^{i}_{k}-(d^{j}-d^{i})\big),\label{eq5}
\end{align}
where $a^{ij}_{k}$ is the coupling coefficient at $t_{k}$. It should be pointed out that the entry $a^{ij}$ of the weighted adjacency matrix $\mathcal{A}_{k}$ allows the existence of negative edge weights. When $a^{ij}_{k}<0$, the dynamics between the agents $i$ and $j$ is competitive. When $a^{ij}_{k}>0$, the dynamics between the agents $i$ and $j$ is cooperative. Moreover, in what follows, the coupling coefficient $a^{ij}_{k}$ is assumed to satisfy: $a^{ij}_{k}\in[\alpha,1)$ and $\sum_{j\in\mathcal{N}^{i}_{k}}a^{ij}_{k}\in[\alpha,1-\alpha]$ for all $a^{ij}_{k}>0$, and $\sum_{j\in\mathcal{N}^{i}_{k}}a^{ij}_{k}\in[-\beta,0)$ for all $a^{ij}_{k}<0$, where $\alpha\in(0,1/2]$ and $\beta\in[0,\alpha/2)$.

Let the position error vector $\varepsilon_{k}^{i}=[\varepsilon_{k}^{i1},\varepsilon_{k}^{i2}]^{T}$ at $t_{k}$ be $\varepsilon_{k}^{i}=p_{k}^{i}-d^{i}$. Combining the dynamics in (\ref{eq1}) and the control protocol in (\ref{eq5}), the position error system is derived by
\begin{align}
\varepsilon^{i}_{k+1}=\varepsilon^{i}_{k}-h_{k}(R^{i}_{k})^{T}\sum_{j=1}^{N} \ell^{ij}_{k}\varepsilon_{k}^{j},\label{eq6}
\end{align}
where $\ell^{ij}_{k}$ is the $ij$-entry of the Laplacian matrix $\mathcal{L}_{k}$.

Due to $\mathcal{L}_{k}\textbf{1}_{N}=\textbf{0}_{N}$, one has $(I_{N}-\mathcal{L}_{k})(I_{N}-\textbf{1}_{N}\textbf{1}^{T}_{N}/N)=I_{N}-\mathcal{L}_{k}-\textbf{1}_{N}\textbf{1}^{T}_{N}/N$. The position error system in (\ref{eq6}) can be rewritten by the vector form
\begin{align}
\varepsilon_{k+1}= (I_{2N}-h_{k}\mathcal{R}_{k}^{T}\Lambda)\varepsilon_{k}+h_{k}\mathcal{R}_{k}^{T}\mathcal{S}_{k}\Lambda\varepsilon_{k},\label{eq7}
\end{align}
where $\varepsilon_{k}=[(\varepsilon^{1}_{k})^{T},\ldots,(\varepsilon^{N}_{k})^{T}]^{T}$, $\Lambda=(I_{N}-\textbf{1}_{N}\textbf{1}^{T}_{N}/N)\otimes I_{2}$, $\mathcal{R}_{k}=\text{blkdiag}\{R^{1}_{k},\ldots,R^{N}_{k}\}$, $\mathcal{S}_{k}=S_{k}\otimes I_{2}$, and $S_{k}=[s^{ij}_{k}]_{N\times N}=I_{N}-\mathcal{L}_{k}\in\mathbb{S}_{\beta}$. The reason introduced the matrix $\Lambda$ into the system in (\ref{eq7}) will be shown in the next section.

\begin{remark}\label{remark1}
As mentioned in \cite{Oh2015}, the normal displacement-based formation control usually has an assumption that each agent can measure the relative positions from neighbors with respect to $\mathfrak{F}^{i}$ due to the limited sensing ability. These relative positions are equivalent to their counterparts with respect to $\mathfrak{F}^{w}$, if $\mathfrak{F}^{w}$ and $\mathfrak{F}^{i}$ are aligned. It further means that the agents require to know the directions of $\mathfrak{F}^{w}$ and $\mathfrak{F}^{i}$. When the direction between $\mathfrak{F}^{w}$ and $\mathfrak{F}^{i}$ is misaligned, a common way is to design a control protocol for estimating and aligning the direction. Although some existing results \cite{Meng2016,Li2019} show that the alignment may be not necessary, there are still some certain restrictions on the rotation angle. Different from the works in \cite{Tran,Lee,Li2019,Meng2016}, this paper considers a class of displacement-based formation control problems, where neither the direction alignments \cite{Tran,Lee} nor the restrictions \cite{Li2019,Meng2016} are required.
\end{remark}

\begin{remark}\label{remark2}
The coupling coefficients express the interactions among the agents in the formation control problems. To realize the desired formation shape, the cooperative interactions are usually necessary to be characterized by the nonnegative coupling coefficients. However, the interactions have not only cooperative properties, but also competitive properties that are represented by the negative coupling coefficients. Actually, the negative coupling coefficients may occur, e.g., some special auxiliary or virtual systems in the theoretical analysis \cite{Altafini,Liu2020}. In this paper, the general coupling coefficients are introduced, which include both nonnegative and negative edge weights. Compared with the works only considered nonnegative coupling coefficients \cite{Dong2019,Zhang2021}, the general coupling coefficients can be endowed with more potential applications.
\end{remark}

\begin{remark}\label{remark4}
Recently, the cooperative and competitive interactions have been considered into the separation problems of formation control \cite{Yan2022,Li2022}. In these works \cite{Yan2022,Li2022}, the agents with cooperative interactions can approach to the desired formation shape, while the agents with competitive interactions can only guarantee the desired relative positions. This further means that the formation shape has two groups separated away from each other under the time-invariant interactions \cite{Altafini}. In contrast to these works \cite{Yan2022,Li2022}, this paper considers the case that the existence of cooperative and competitive interactions does not affect the desired formation shape, where a time-varying signed graph is used to characterize the interactions.
\end{remark}

\section{Formation control with unknown direction information and general topology}\label{S3}

The formation control problem with unknown directions and general coupling coefficients in (\ref{eq7}) is studied in this section. First, the properties of several polytopes are considered, which contribute to derive the main result.

\subsection{Technical analysis of several polytopes}\label{S31}

In this subsection, some polytopes are constructed, and their properties are discussed on a certain hyperplane. According to the contracting mapping principle, these polytopes are then indicated to converge under the mappings of the matrices $h_{k}\mathcal{R}_{k}^{T}\Lambda$ and $h_{k} \mathcal{R}_{k}^{T}\mathcal{S}_{k}\Lambda$, respectively.

Now, construct an $\mathcal{H}$-polytope
\begin{align}
\mathcal{P}_{N}=\bigcap_{i=1}^{2^N}\{x\in\mathbb{R}^{N}|\langle a^{i},x\rangle\in[-1,1]\},\label{eq8}
\end{align}
where $a^{i}\in\mathbb{R}^{N}$ is a vector whose entry $a^{ip}$ satisfies $a^{ip}\in\{\pm1\}$ for $p=\{1,\ldots,N\}$, and $a^{i}\neq a^{j}$ for $i\neq j$. Apparently, the polytope $\mathcal{P}_{N}$ is equivalent to the form
\begin{align}
\mathcal{P}_{N}=\{x\in\mathbb{R}^{N}|\|x\|_{1}\leq 1\},\label{eq9}
\end{align}
which is called \emph{crosspolytope}. Then, the following property is well known \cite{Ziegler}, and its proof is omitted.

\begin{proposition}\label{proposition1}
$\text{ext}\mathcal{P}_{N}=\{\pm e_{N}^{i}\}_{i=1}^{N}$.
\end{proposition}

$\mathcal{H}=\{x\in\mathbb{R}^{N}| \langle\textbf{1}_{N},x\rangle  =0\}$ is utilized to denote as a hyperplane, and an intersection associate with $\mathcal{P}_{N}$ and $\mathcal{H}$ is defined by $
\mathcal{P}'_{N}=\mathcal{P}_{N}\cap\mathcal{H}$. The following lemma solves $\text{ext}\mathcal{P}'_{N}$.

\begin{lemma}\label{lemma1}
$\text{ext}\mathcal{P}'_{N}=\mathcal{V}$, where $\mathcal{V}=\{\widetilde{e}^{i}_{N}\}_{i=1}^{\widetilde{N}}$, $\widetilde{N}$ is an even number, $\widetilde{e}^{i}_{N}=(e^{p}_{N}-e^{q}_{N})/2$ ($p\neq q$ and $p,q\in\{1,\ldots,N\}$), $\widetilde{e}^{i}_{N}\neq\widetilde{e}^{j}_{N}$ for $i\neq j$, and $\widetilde{e}^{i}_{N}=-\widetilde{e}^{\widetilde{N}+1-i}_{N}$ by arranging the order of $\widetilde{e}^{i}_{N}$.
\end{lemma}
\begin{IEEEproof}
See the Appendix \ref{appendix1}.
\end{IEEEproof}

Lemma \ref{lemma1} describes the equivalence of $\mathcal{P}'_{N}$ and $\mathcal{P}_{\mathcal{V}}$, whose extreme point sets are same, i.e., $\mathcal{V}$. However, Lemma \ref{lemma1} cannot be applied to analyzing $\mathcal{P}_\mathcal{V}$ under the mappings of the matrices $h_{k}\mathcal{R}_{k}^{T}\Lambda$ and $h_{k}\mathcal{R}_{k}^{T}\mathcal{S}_{k}\Lambda$, respectively. Similarly, construct a polytope $\overline{\mathcal{P}}_{2N}=\mathcal{P}_{2N}\cap\overline{\mathcal{H}}$, where $\overline{\mathcal{H}}=\{x\in\mathbb{R}^{2N}| \langle\textbf{1}_{N}\otimes e_{2}^{i},x\rangle  =0\}_{i=1}^{2}$. Now, the result in Lemma \ref{lemma1} is extended into a higher dimension space as follows.

\begin{lemma}\label{lemma2}
$\text{ext}\overline{\mathcal{P}}_{2N}=\overline{\mathcal{V}}$, where $\overline{\mathcal{V}}=\{\overline{e}_{N,j}^{i}| \overline{e}_{N,j}^{i}=\widetilde{e}_{N}^{i}\otimes e^{j}_{2}\}_{i=1}^{\widetilde{N}}$ for $j=1,2$, and $\widetilde{e}_{N}^{i}$ is defined in Lemma \ref{lemma1}.
\end{lemma}
\begin{IEEEproof}
See the Appendix \ref{appendix2}.
\end{IEEEproof}

Then, consider $\mathcal{S}_{k}^{T}$ to map $\mathcal{P}_{\overline{\mathcal{V}}}$ and $\mathcal{P}_{2N}$ on $\overline{\mathcal{H}}$, respectively. The following results are derived.

\begin{lemma} \label{lemma3}
For any $k\in\mathbb{N}$, $\mathcal{S}_{k}^{T}\mathcal{P}_{\overline{\mathcal{V}}}\subset\mathcal{P}_{\overline{\mathcal{V}}}$, if $\mathcal{G}_{k}$ is neighbor shared.
\end{lemma}
\begin{IEEEproof}
See the Appendix \ref{appendix3}.
\end{IEEEproof}

\begin{lemma} \label{lemma4}
For any $k\in\mathbb{N}$, $\mathcal{S}_{k}^{T}\mathcal{P}_{\overline{\mathcal{V}}}=\mathcal{S}_{k}^{T}\mathcal{P}_{2N}\cap\overline{\mathcal{H}}$, if $\mathcal{G}_{k}$ is neighbor shared.
\end{lemma}
\begin{IEEEproof}
See the Appendix \ref{appendix4}.
\end{IEEEproof}

Lemma \ref{lemma3} derives a result that the mapping $\mathcal{S}_{k}^{T}\mathcal{P}_{\overline{\mathcal{V}}}$ is a contractive set related to $\mathcal{P}_{\overline{\mathcal{V}}}$. Lemma \ref{lemma4} further implies that this contractility can be valid for the mapping $\mathcal{S}_{k}^{T}\mathcal{P}_{2N}$ on $\overline{\mathcal{H}}$. However, when $\mathcal{R}_{k}\neq I_{2N}$, the similar results are impossible for the mapping $\mathcal{R}_{k}\mathcal{P}_{\overline{\mathcal{V}}}$ due to $\langle\textbf{1}_{N}\otimes e_{2}^{p},\mathcal{R}_{k}\overline{e}_{N,j}^{i}\rangle  \neq0$. To approximate the action region of $\mathcal{R}_{k}\mathcal{P}_{\overline{\mathcal{V}}}$, denote the polar of $\mathcal{P}_{2N}$ as $\mathcal{P}_{2N}^{\circ}$. The following propositions reveal the relationship among $\mathcal{P}_{\overline{\mathcal{V}}}$, $ \mathcal{P}_{2N}$ and $\mathcal{P}_{2N}^{\circ}$.

\begin{proposition}\label{proposition2}
For any $k\in\mathbb{N}$, $\mathcal{R}_{k}\mathcal{P}_{\overline{\mathcal{V}}}\subseteq \mathcal{R}_{k}\mathcal{P}_{2N}\subseteq \mathcal{P}_{2N}^{\circ}$.
\end{proposition}
\begin{IEEEproof}
See the Appendix \ref{appendix5}.
\end{IEEEproof}

\begin{proposition}\label{proposition3}
For any $k\in\mathbb{N}$, $h_{k}\mathcal{P}_{2N}^{\circ}\subseteq \mathcal{P}_{2N}$, if $h_{k}\in(0,1/2N]$.
\end{proposition}
\begin{IEEEproof}
See the Appendix \ref{appendix6}.
\end{IEEEproof}

Propositions \ref{proposition2} and \ref{proposition3} express a way that should be able to construct a polytope based on $\mathcal{P}_{2N}^{\circ}$, such that the mapping $\mathcal{R}_{k}\mathcal{P}_{\overline{\mathcal{V}}}$ could be treated as the similar results of Lemmas \ref{lemma3} and \ref{lemma4}. In order to do this, the matrix $\Lambda$ in (\ref{eq7}) is considered as follows.

\begin{proposition} \label{proposition4}
For any $k\in\mathbb{N}$, $\Lambda\mathcal{R}_{k}\mathcal{P}_{\overline{\mathcal{V}}}\subseteq\mathcal{P}_{2N}^{\circ}\cap\overline{\mathcal{H}}$.
\end{proposition}
\begin{IEEEproof}
See the Appendix \ref{appendix7}.
\end{IEEEproof}

In Proposition \ref{proposition4}, $\Lambda$ is actually an orthogonal projection matrix. Therefore, Proposition \ref{proposition4} indicates a geometrical intuition that $\Lambda$ projects the mapping $\mathcal{R}_{k}\mathcal{P}_{\overline{\mathcal{V}}}$ on $\overline{\mathcal{H}}$. This is why the matrix $\Lambda$ is introduced into the system in (\ref{eq7}).

Now, from the above results, the mappings $h_{k}\Lambda\mathcal{R}_{k}\mathcal{P}_{\overline{\mathcal{V}}}$ and $h_{k} \Lambda\mathcal{S}_{k}^{T}\mathcal{R}_{k}\mathcal{P}_{\overline{\mathcal{V}}}$ are considered by the following lemma.
\begin{lemma} \label{lemma5}
For any $k\in\mathbb{N}$, $h_{k}\Lambda\mathcal{R}_{k} \mathcal{P}_{\overline{\mathcal{V}}}\subseteq\mathcal{P}_{\overline{\mathcal{V}}}$ and $h_{k} \Lambda\mathcal{S}_{k}^{T}\mathcal{R}_{k}$ $\mathcal{P}_{\overline{\mathcal{V}}}\subseteq\mathcal{P}_{\overline{\mathcal{V}}}$, if all the conditions in Lemma \ref{lemma3} and Proposition \ref{proposition3} hold.
\end{lemma}
\begin{IEEEproof}
Obviously, $h_{k} \Lambda\mathcal{S}_{k}^{T}\mathcal{R}_{k}\mathcal{P}_{\overline{\mathcal{V}}}\subseteq h_{k} \Lambda\mathcal{S}_{k}^{T}\mathcal{R}_{k}\mathcal{P}_{\overline{\mathcal{V}}}\cap\overline{\mathcal{H}}$, and $\Lambda\mathcal{P}_{\overline{\mathcal{V}}}\subseteq \mathcal{P}_{\overline{\mathcal{V}}}$. Then, the proof directly follows from Lemmas \ref{lemma3}-\ref{lemma4} and Propositions \ref{proposition2}-\ref{proposition4}, which is omitted.
\end{IEEEproof}

\begin{remark}
The aim of this subsection is to construct some polytopes, such that the difference of any pair of rows of $ \mathcal{R}_{k}$ and $\mathcal{S}_{k}^{T}$ can be analyzed. This difference plays an important factor for dealing with the cooperative control problems. Lemma \ref{lemma1} can be utilized to describe this difference, whose counterpart in the higher dimension space is extended by Lemma \ref{lemma2}. According to the contractive set method, Lemmas \ref{lemma3} and \ref{lemma4} indicate that the mapping $\mathcal{S}_{k}^{T}\mathcal{P}_{\overline{\mathcal{V}}}$ is contractive, and is further equivalent to the inclined projection of $\mathcal{S}_{k}^{T}\mathcal{P}_{2N}$ on $\overline{\mathcal{H}}$. However, the mapping $\mathcal{R}_{k} \mathcal{P}_{\overline{\mathcal{V}}}$ is intuitively impossible to maintain the similar contractive property. In this case, several properties of $\mathcal{P}_{2N}^{\circ}$ are considered by Propositions \ref{proposition2}-\ref{proposition4}. As a result, Lemma \ref{lemma5} solves the contractive property of the mappings $h_{k}\Lambda\mathcal{R}_{k} \mathcal{P}_{\overline{\mathcal{V}}}$ and $h_{k} \Lambda\mathcal{S}_{k}^{T}\mathcal{R}_{k}\mathcal{P}_{\overline{\mathcal{V}}}$. It should be pointed out that $h_{k}\mathcal{P}_{2N}^{\circ}$ can be regraded as a contraction mapping, and $\Lambda$ is applied to projecting $\mathcal{S}_{k}^{T}\mathcal{R}_{k}\mathcal{P}_{\overline{\mathcal{V}}}$ on $\overline{\mathcal{H}}$.
\end{remark}

\subsection{The convergence of formation control}\label{S32}

According to the discussions in the above subsection, the convergence problem of the system in (\ref{eq7}) will be considered. Before presenting the main result, let $\mathcal{E}_{N}=E_{N}\otimes I_{2} $, where $E_{N}=[\widetilde{e}_{N}^{1},\ldots,\widetilde{e}_{N}^{\#\mathcal{V}}]$. Then, one obtains the following result.

\begin{lemma}\label{lemma6}
Assume that all the conditions in Lemma \ref{lemma5} hold. If there exists a nonnegative matrix $\mathcal{M}_{k}$, then the following condition holds
\begin{align*}
&\Big(I_{2N}-h_{k}(\mathcal{L}_{k}^{T}\otimes I_{2})\mathcal{R}_{k} \Big) \mathcal{E}_{N}=\mathcal{E}_{N}\mathcal{M}_{k},
\end{align*}
where $\|\mathcal{M}_{k}\|_{1}\leq   1-(1-\eta)h_{k}$, and $\eta=1-\alpha+2\beta<1$.
\end{lemma}
\begin{IEEEproof}
Due to $S_{k}=I_{N}-\mathcal{L}_{k}$ and $\mathcal{L}_{k}\textbf{1}_{N}=\textbf{0}_{N}$, one has
\begin{align*}
 I_{2N}-h_{k} (\mathcal{L}_{k}^{T}\otimes I_{2})\mathcal{R}_{k}=& I_{2N}-h_{k}\Lambda(\mathcal{L}_{k}^{T}\otimes I_{2})\mathcal{R}_{k}\\
=&  I_{2N}-h_{k}\Lambda\mathcal{R}_{k} +h_{k}\Lambda\mathcal{S}_{k}^{T}\mathcal{R}_{k}.
\end{align*}
Then, one can use $x\in\mathcal{P}_{\overline{\mathcal{V}}}$ in place of $y\in\mathcal{P}_{2N}$ in (\ref{eqp41}). Thus, the proof follows from the discussions of Lemmas \ref{lemma3}-\ref{lemma5} and Propositions \ref{proposition2}-\ref{proposition4}.
\end{IEEEproof}

Based on Lemma \ref{lemma6}, for any pair of distinct agents $i,j$, the relative position vector $p^{j}_{k}-p^{i}_{k}$ is proven to tend to the desired relative position vector $d^{j}-d^{i}$.

\begin{theorem}\label{theorem1}
Assume that all the conditions in  Lemma \ref{lemma6} hold. If there exists a positive scalar $\delta\in(0,1)$, such that the sampling interval $h_{k}\in[h^{m},h^{M}]$ is solvable, then the formation control problem in (\ref{eq7}) is achieved under the convergence rate $-\ln(1-\delta)/h^{M}$, where  $h^{M}\leq1/2N$, $h^{m}\geq\delta/(1-\eta)$, $h^{m}$ and $h^{M}$ stand for the upper and lower bounds of $h_{k}$, respectively.
\end{theorem}

\begin{IEEEproof}
From the position error system in (\ref{eq7}), one has
\begin{align}
 \varepsilon_{k+1}=&  (I_{2N}-h_{k}\mathcal{R}_{k}^{T}\Lambda +h_{k}\mathcal{R}_{k}^{T}\mathcal{S}_{k}\Lambda )\varepsilon_{k}\nonumber\\
=& \prod_{j=0}^{k}\Big( I_{2N}-h_{j}\mathcal{R}_{j}^{T}(\mathcal{L}_{j}\otimes I_{2})\Lambda  \Big)\varepsilon_{0}.
\label{eq10}
\end{align}
where $\prod$ means the left matrix products.

According to Lemma \ref{lemma6}, one gets
\begin{align}
\|\mathcal{E}&_{N}^{T}\varepsilon_{k+1}\|_{\infty}\nonumber\\
\leq&\Big\|\mathcal{E}_{N}^{T}\prod_{j=0}^{k}\Big( I_{2N}-h_{j}\mathcal{R}_{j}^{T}(\mathcal{L}_{j}\otimes I_{2})\Lambda  \Big)\Big\|_{\infty}
 \| \varepsilon_{0 }\|_{\infty}\nonumber\\
=&\Big\| \prod_{j=0}^{k} \Big( I_{2N}-h_{j}\Lambda(\mathcal{L}_{j}^{T}\otimes I_{2})\mathcal{R}_{j}\Big)\mathcal{E}_{N}\Big\|_{1} \| \varepsilon_{0 }\|_{\infty}\nonumber\\
\leq&\prod_{j=0}^{k}\Big(  1-(1-\eta)h_{j}\Big)\|\mathcal{E}_{N} \|_{1}\|\varepsilon_{0}\|_{\infty},\label{eq11}
\end{align}
where $\eta$ is given in Lemma \ref{lemma6}.

Due to $h^{m}\geq\delta/(1-\eta)$ and $h^{M}\leq1/2N$, one has $1-(1-\eta) h_{k}\leq1-\delta$, if $\delta\in(0,\frac{1-\eta}{2N}]$. In case $t_{k+1}-t_{k}\leq h^{M}$, then $t_{k+1}\leq t_{0}+(k+1)h^{M}$. From (\ref{eq11}), it yields
\begin{align*}
\|\mathcal{E}_{N}^{T}\varepsilon_{k+1}\|_{\infty} \leq&(1-\delta)^{ k+1 }\|\mathcal{E}_{N} \|_{1}\|\varepsilon_{0}\|_{\infty}\\
\leq&\exp\Big\{ \frac{\ln(1-\delta)}{h^{M}}(t_{k+1}-t_{0})\Big\}\|\mathcal{E}_{N} \|_{1}\|\varepsilon_{0}\|_{\infty}.
\end{align*}
Therefore, one gets
\begin{align*}
\lim_{k\rightarrow\infty}\|\mathcal{E}_{N}^{T}\varepsilon_{k}\|_{\infty}=0.
\end{align*}

Recalling the structure of $\mathcal{E}_{N}=E_{N}\otimes I_{2} $, the $i$-column of $E_{N}$ is $ \widetilde{e}_{N}^{i} $, which further means that
\begin{align*}
\lim_{k\rightarrow\infty}\|p^{p}_{k}-p^{q}_{k}-d^{p}+d^{q}\|_{\infty}=0.
\end{align*}
It completes the proof.
\end{IEEEproof}

\begin{remark}
In many displacement-based formation control problems \cite{Marina,Dong2019,Zhang2021}, the Lyapunov function usually plays an important part. However, in some cases \cite{Olshevsky,Matei}, e.g., the dynamic interactions or the nonexistence of eigenvalues of system matrices, it may be hard to construct a common Lyapunov function whose convergence may not also be guaranteed. This paper introduces the unknown matrix $\mathcal{R}_{k}^{T}$ and the general stochastic matrix $\mathcal{S}_{k}$ in (\ref{eq7}), such that the construction and analysis of the Lyapunov function are inevitably more complicated. In this case, several polytopes are applied to considering the convergence of $\mathcal{S}_{k}$. Then, a larger convex region induced by the polar is to enclose the action of $\mathcal{R}_{k}$. According to the contracting mapping principle and the obtained results, the convergence problem of the system in (\ref{eq7}) is finally solved. Different from the Lyapunov function method \cite{Marina,Dong2019,Zhang2021}, the approach in this paper relies on geometrical intuition, which means that some restrictions of the Lyapunov function are removed.
\end{remark}

\section{Example}
An simulation of the formation control with the unknown direction information and the general topology is provided in this section. The system has $4$ agents moving in $\mathbb{R}^{2}$ ($\mathbb{R}^{2}$ is assume to be $\mathfrak{F}^{w}$).

\begin{figure}[htbp!]
  \centering
  \includegraphics[width=3in]{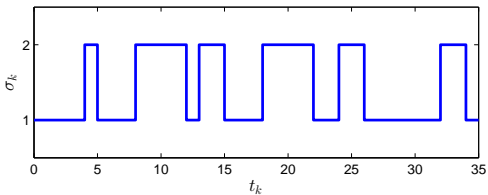}\\
  \caption{The switching signal $\sigma_{k}=\{1,2\}$.}\label{fig:1}
\end{figure}
\begin{figure}[htbp]
\centering
\subfigure[The graph $\mathcal{G}_{k}$ for $\sigma_{k}=1$.]
{\includegraphics[height=4cm,width=4cm]{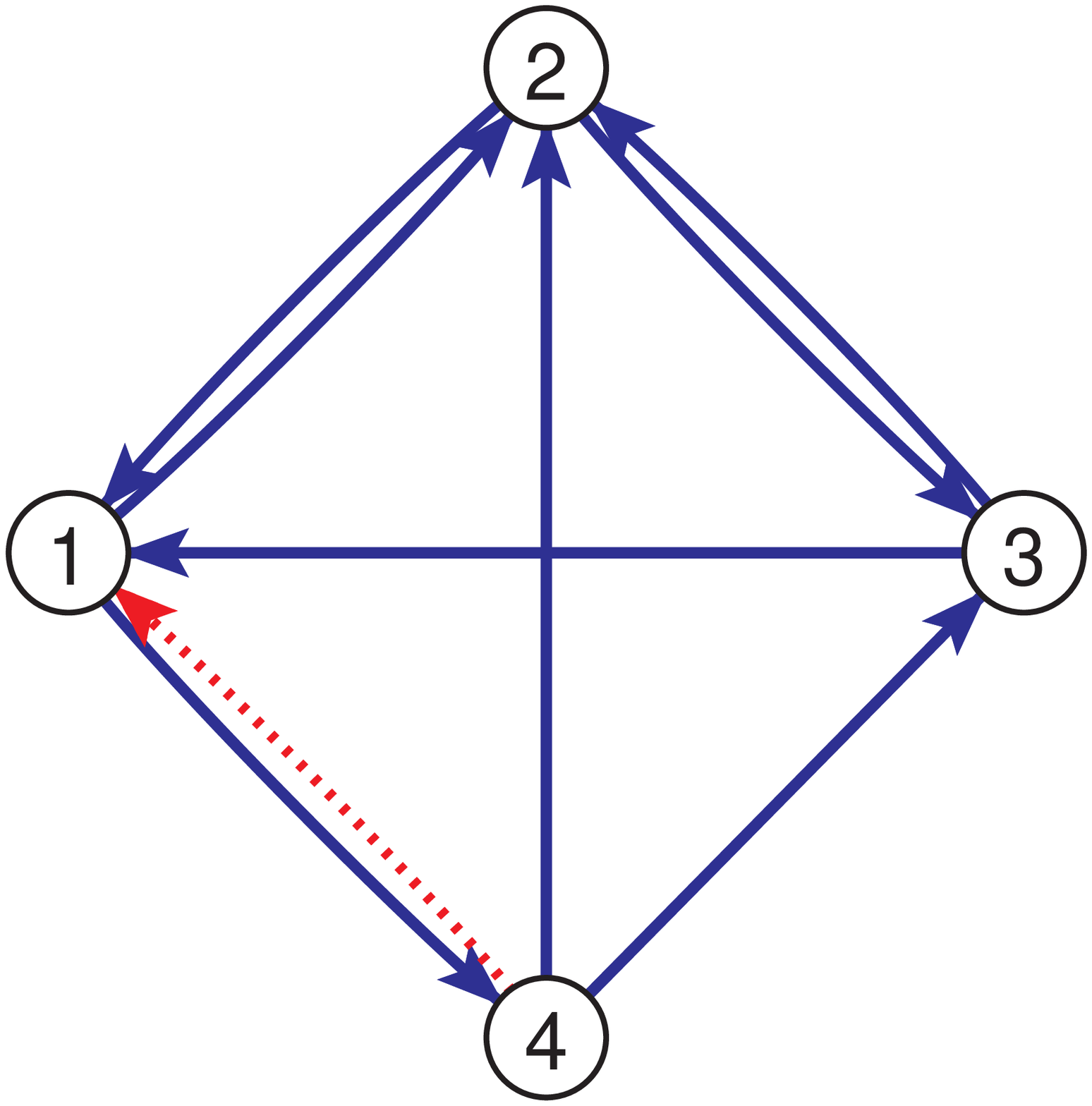}}
\subfigure[The graph $\mathcal{G}_{k}$ for $\sigma_{k}=2$.]
{\includegraphics[height=4cm,width=4cm]{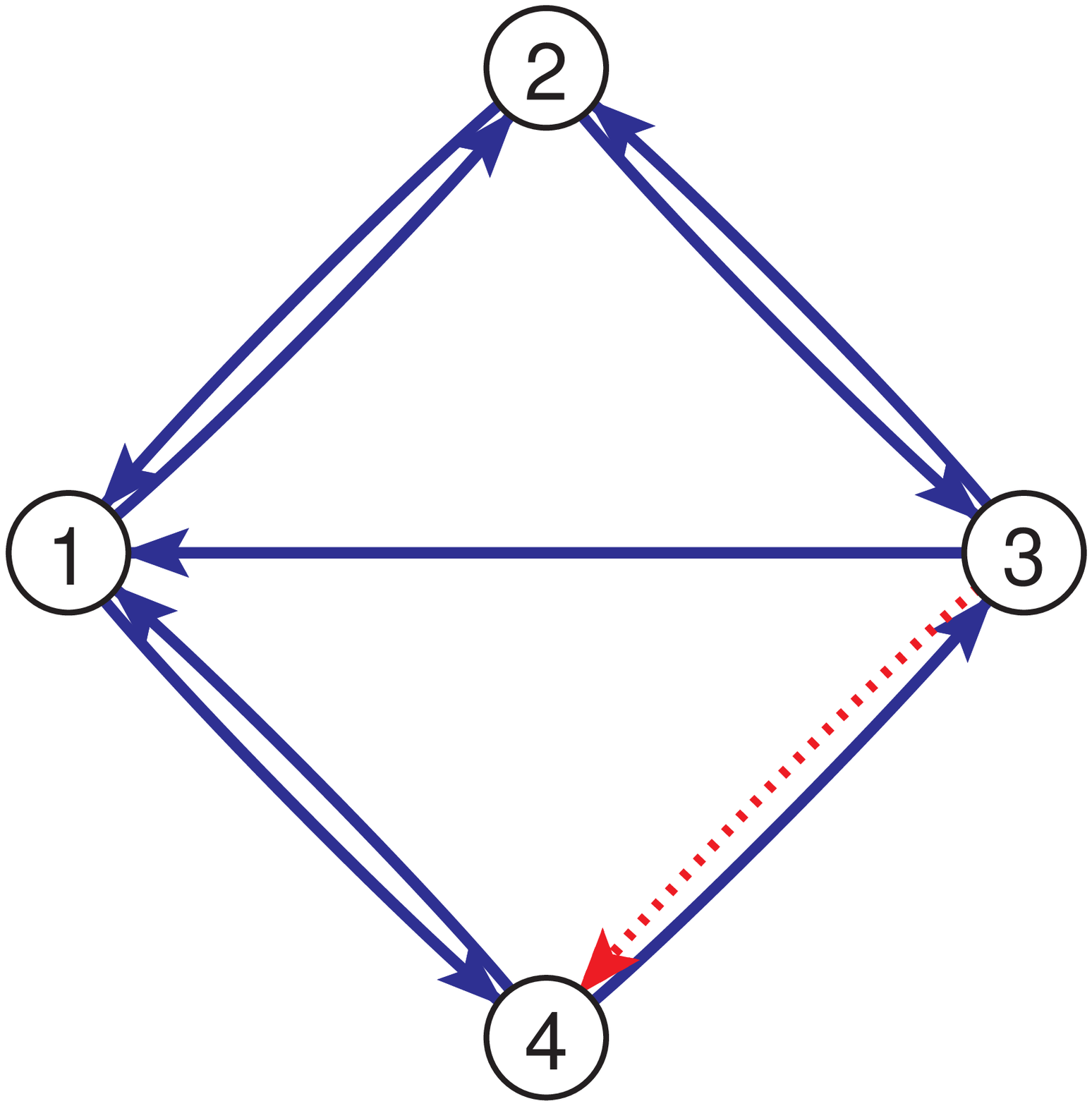}}
\caption{The signed directed graph $\mathcal{G}_{k}$, where the blue arcs are the positive edges with weights $0.2 $, and the red dashed arcs are the negative edges with weights $-0.08 $. }
\label{fig:2}
\end{figure}

For the formation control protocol in (\ref{eq5}), the parameters are designed in the following. The absolute desired position vector $d$ with respect to $\mathfrak{F}^{w}$ is given as $d^{1}=[0,0]^{T}$, $d^{2}=[5,0]^{T}$, $d^{3}=[5,5]^{T}$, and $d^{4}=[0,5]^{T}$, i.e., the desired formation shape is a square with length $5$. The position vector $p_{k}$ with respect to $\mathfrak{F}^{w}$ is initialized by $p^{1}_{0}=[-6,2]^{T}$, $p^{2}_{0}=[9,-10]^{T}$, $p^{3}_{0}=[3,0]^{T}$ and $p^{4}_{0}=[0,16]^{T}$. Fig. \ref{fig:1} shows that a switching signal $\sigma_{k}=\{1,2\}$ determines the time-varying directed graph $\mathcal{G}_{k}$. Fig. \ref{fig:2} draws the signed graph $\mathcal{G}_{k}$, where the weights of the positive edges (the blue arcs) and the negative edges (the red dashed arcs) are all equal to $0.2$ and $-0.08$, respectively.

\begin{figure}[htbp!]
  \centering
  \includegraphics[width=3in]{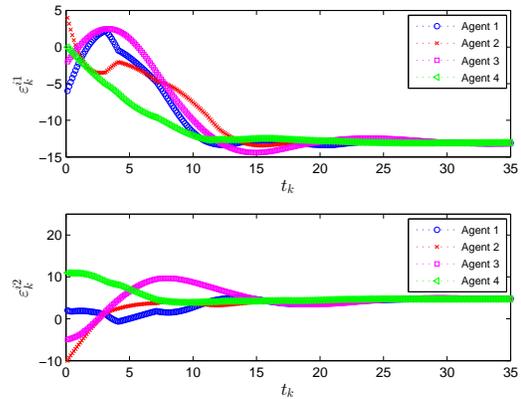}\\
  \caption{The position error trajectories $\varepsilon_{k}^{i}=p_{k}^{i}-d^{i}$, where $h_{k}=0.125$ for all $k\in\mathbb{N}$.}\label{fig:3}
\end{figure}

\begin{figure}[htbp!]
  \centering
  \includegraphics[width=3in]{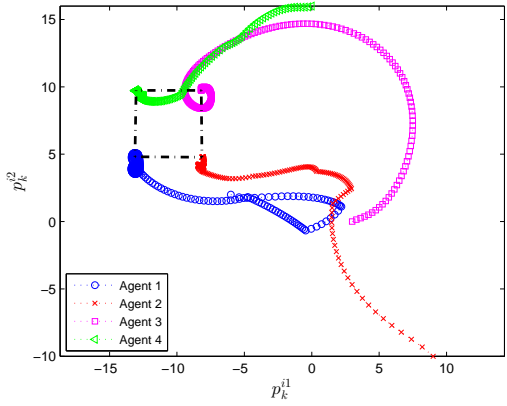}\\
  \caption{The desired formation shape and position trajectories of multi-agent system, where the initial value of $p_{k}$ are $p^{1}_{0}=[-6,2]^{T}$, $p^{2}_{0}=[9,-10]^{T}$, $p^{3}_{0}=[3,0]^{T}$, $p^{4}_{0}=[0,16]^{T}$, and the black dashed lines are the desired formation shape,.}\label{fig:4}
\end{figure}

Although, for each agent $i$ and $k\in\mathbb{N}$, the rotation angle $\theta_{k}^{i}$ can be unknown, $\theta_{k}^{i}$ is assumed to be the following function in order to facilitate the simulation with respect to $\mathfrak{F}^{w}$: $\theta_{k}^{i}=\frac{\pi t_{k}}{2i(1+|t_{k}|)}$ for $i=1,2$, and $\theta_{k}^{i}=-\frac{i\pi  t_{k}}{8(1+|t_{k}|)}$ for $i=3,4$. From Theorem \ref{theorem1}, one has $h_{k}\in[ 100\delta, 0.125]$ and $\eta=1-\alpha+2\beta=0.96$, which further results in $\delta\in(0, 0.12]$. In this case, the desired formation shape is realized, since $\delta$ is solvable. One calculates the convergence rate to satisfy: $-\ln(1-\delta)/h^{M}\in(0,0.128]$. For all $i=1,\ldots,4$ and $k\in\mathbb{N}$, the position errors $\varepsilon_{k}^{i}=p_{k}^{i}-d^{i}$ are expressed by selecting $h_{k}=h^{M}=0.125$ in Fig. \ref{fig:3}. Based on the same selection of $h_{k}$, the desired formation shape and the position trajectories of agents are drawn in Fig. \ref{fig:4}.

\section{Conclusion}
This paper investigates a kind of displacement-based formation control problem. Each agent $i$ is assumed to gather the measurements from the neighbors with respect to $\mathfrak{F}^{i}$, but cannot directly know the direction of $\mathfrak{F}^{i}$. Meanwhile, each agent $i$ is also assumed to communicate with the neighbors by utilizing the general coupling coefficients that describe the cooperative and competitive interactions for the whole system. For the sake of dealing with the above assumptions, several convex polytopes have been introduced. The properties of the polytopes are then analyzed  in view of some mathematical techniques, which further reveals that the mappings of the unknown matrix $\mathcal{R}_{k}^{T}$ and the general stochastic matrix $\mathcal{S}_{k}$ on these polytopes are contractive. By utilizing the derived results, the desired formation shape is finally guaranteed, while the cooperative and competitive interactions coexist under the measurements from unknown direction of each agent. Subsequently, an example is simulated for demonstrating the final results. Further research interests involve the displacement-based formation control on $\mathrm{SO}(3)$.

\section{Appendices}

\begin{appendices}
Some technical results of Subsection \ref{S31} are proven in these appendices.
\section{Proof of Lemma \ref{lemma1}}\label{appendix1}
\begin{IEEEproof}
Solving $\text{ext}\mathcal{P}'_{N}$ is equivalent to prove that $\mathcal{P}'_{N}$ is equivalent to $\mathcal{P}_\mathcal{V}$, and any point in $ \mathcal{V}$ cannot be expressed by a convex combination of at least two distinct points in $\mathcal{P}_\mathcal{V}$. Then, the proof has two steps.

\emph{Step I}: To prove $\mathcal{P}_\mathcal{V}=\mathcal{P}'_{N}$.

For any $x\in\mathcal{P}_{\mathcal{V}}$, $x$ has a form of a convex combination of the entries in $\mathcal{V}$, i.e., $x=\sum_{i=1}^{\widetilde{N}}\lambda^{i}\widetilde{e}^{i}_{N}\in\mathcal{P}_{\mathcal{V}}$, where $\sum_{i=1}^{\widetilde{N}}\lambda^{i}=1$ and $\lambda^{i}\in[0,1]$. Construct a matrix $A_{N}=[a^{1},\ldots,a^{2^{N}}]^{T}\in \mathbb{R}^{2^{N}\times N}$, where $a^{i}$ is given in (\ref{eq8}). The polytope $\mathcal{P}_{N}$ can be rewritten by the form $\mathcal{P}_{N}=\{x\in\mathbb{R}^{N}|A_{N}x\leq \textbf{1}_{2^N}\}$. Then, one gets
\begin{align*}
A_{N}x=\sum_{i=1}^{\widetilde{N}}\lambda^{i}A_{N}\widetilde{e}^{i}_{N}\leq\sum_{i=1}^{\widetilde{N}}\lambda^{i}\textbf{1}_{2^N}=\textbf{1}_{2^N},
\end{align*}
where $A_{N}\widetilde{e}^{i}_{N}\leq\textbf{1}_{2^N}$ holds due to $\langle a^{j},\widetilde{e}^{i}_{N}\rangle\in[-1,1]$ for all $a^{j}$ in (\ref{eq8}), $j=1,\ldots,2^N$ and $i=1,\ldots,\widetilde{N}$. In this case, $x\in\mathcal{P}_{N}$ holds. Moreover, note that $\langle\textbf{1}_{N},\widetilde{e}^{i}_{N}\rangle  =0$, it yields $x\in\mathcal{H}$. Therefore, $\mathcal{P}_{\mathcal{V}}\subseteq\mathcal{P}'_{N}$ holds.

For the converse, the proof relies on a geometrical intuition that suggests to solve an inclined projection of $\mathcal{P}_{N}$ on $\mathcal{H}$. Then, the mathematical induction is utilized.

When $N=2$, there is no difference between the inclined projection and the orthogonal projection of $\mathcal{P}_{2}$ onto $\mathcal{H}$. Thus, the Fourier-Motzkin elimination can be utilized by eliminating the variable along the axis $\textbf{1}_{2}$.

Note that all the solutions of $A_{2}x\leq\textbf{1}_{4}$ are $\mathcal{P}_{2}$, where $A_{2}= \big[\begin{smallmatrix}
1&  1&  -1&  -1\\
1& -1&  -1&  1
\end{smallmatrix}\big]^{T}$. Let $\overline{x}=[\overline{x}^{1},\overline{x}^{2}]^{T}\in\text{span}\{e_{2}^{i}\}_{i=1}^{2}$ and $\overline{y}=[\overline{y}^{1},\overline{y}^{2}]^{T}\in\text{span}\{\textbf{1}_{2}/\sqrt{2},(e^{1}_{2}-e^{2}_{2})/\sqrt{2}\}$, one can construct an affine mapping from $\overline{x}$ to $\overline{y}$, such that the following inequality holds for any $\overline{x}\in\mathcal{P}_{2}$
\begin{align*}
\left[\begin{array}{c}
I_{2}\\
-I_{2}\\
\end{array}\right]
\left[\begin{array}{c}
\overline{y}^{1}\\
\overline{y}^{2}\\
\end{array}\right]=
\left[\begin{array}{c}
I_{2}\\
-I_{2}\\
\end{array}\right]
\left[\begin{array}{c}
\overline{x}^{1}+\overline{x}^{2}\\
\overline{x}^{1}-\overline{x}^{2}\\
\end{array}\right]=A_{2}\overline{x}\leq\textbf{1}_{4}.
\end{align*}
By eliminating the variable along the axis $\textbf{1}_{2}$, i.e., $\overline{y}^{1}$, one has $ \overline{y}^{2}\in[-1,1]$,
which further implies $\text{bd}\overline{x}\in\{\pm( e_{2}^{1}-e_{2}^{2})/2\}$ due to the definition of $\mathcal{H}$, where $\text{bd}\overline{x}$ is the boundary of $\overline{x}$.

When $N=K$, assume that the boundary of solutions is in $\mathcal{V}$. Then, when $N=K+1$, $\mathcal{P}_{K+1}$ has the form $A_{K+1}x\leq\textbf{1}_{2^{K+1}}$, which is equivalent to
\begin{align}
\left[\begin{array}{cc}
A_{K}&1\\
A_{K}&-1\\
\end{array}\right]
\left[\begin{array}{c}
\widetilde{x}\\
\widehat{x}
\end{array}\right]\leq
\textbf{1}_{2^{K+1}},\label{eql11}
\end{align}
where $\widetilde{x}\in\mathbb{R}^{K}$, $\widehat{x}\in\mathbb{R}$, and $x=[\widetilde{x}^{T},\widehat{x}]^{T}$.

If $\widehat{x}=0$, it is easy to the result due to the case of $N=K$.

If $\widehat{x}\neq0$, one obtains $\sum_{i=1}^{K}\widetilde{x}^{i}=-\widehat{x}$, since the solutions of $x$ are on $\mathcal{H}$. From (\ref{eql11}), one has $ A_{K}\widetilde{x}\leq\textbf{1}_{2^{K}}/2$. According to Proposition \ref{proposition1}, the boundary of solutions of $\widetilde{x}$ is in $\text{ext}\{x\in\mathbb{R}^{K}|A_{K}x\leq \textbf{1}_{2^K}/2\}$, which means $\text{bd}\widetilde{x}\in\{\pm e_{K}^{i}/2\}_{i=1}^{K}$ and $\text{bd}\widehat{x}\in\{\pm 1/2\} $. In this case, it is not hard to obtain $\text{bd}\big[\begin{smallmatrix}
\widetilde{x}\\
\widehat{x}
\end{smallmatrix}\big]\in\{\pm( e_{K+1}^{i}-e_{K+1}^{K+1})/2\}_{i=1}^{K}$. Therefore, $\mathcal{P}'_{N}\subseteq\mathcal{P}_{\mathcal{V}}$ holds.

\emph{Step II}: To prove $\text{ext}\mathcal{P}_\mathcal{V}=\mathcal{V}$.

If there exists an extreme point $v\notin\mathcal{V}$, then there does not exist any $\lambda^{i}\in[0,1]$, such that $v=\sum_{i=1}^{\widetilde{N}}\lambda^{i}\widetilde{e}_{N}^{i}$, where $\sum_{i=1}^{\widetilde{N}}\lambda^{i}=1$. It is equivalent to that there does not exist a vector $\lambda=[\lambda^{1},\ldots,\lambda^{\widetilde{N}}]^{T}$, such that the following equation holds
\begin{align*}
\left[\begin{array}{c}
\textbf{1}_{\widetilde{N}}^{T}\\
\sum_{i=1}^{\widetilde{N}}\widetilde{e}_{N}^{i}\otimes (e_{\widetilde{N}}^{i})^{T}\\
\end{array}\right]\lambda=
\left[\begin{array}{c}
1\\
v\\
\end{array}\right] .
\end{align*}

From the Frakes Lemma II \cite{Ziegler}, there exists a column vector $\phi$ with appropriate dimensions, such that $\Big\langle \phi, \big[\begin{smallmatrix}
1\\
v
\end{smallmatrix}\big]\Big\rangle<0$ and $ \phi^{T} \Big[\begin{smallmatrix}
\textbf{1}_{\widetilde{N}}^{T}\\
\sum_{i=1}^{\widetilde{N}}\widetilde{e}_{N}^{i}\otimes (e_{\widetilde{N}}^{i})^{T}
\end{smallmatrix}\Big]\geq \textbf{0}_{\widetilde{N}}^{T}$ hold. Let $\phi=[\varphi,-\widetilde{\varphi}^{T}]^{T}$, $\varphi\in\mathbb{R}$ and $\widetilde{\varphi}\in\mathbb{R}^{N}$. One gets $\langle \widetilde{\varphi}, v\rangle>\varphi$ and $\langle \widetilde{\varphi}, \widetilde{e}_{N}^{i}\rangle\leq\varphi$. Then, the proof has two cases.

\emph{Case I}: $\varphi$ is the largest signed distance along the direction $\widetilde{\varphi} $ from the origin of $ \mathcal{P}_\mathcal{V}$ to a hyperplane $\mathcal{H}'$ containing $\widetilde{e}_{N}^{i}$. It is to say that $\mathcal{H}'$ supports $\mathcal{P}_\mathcal{V}$ at $\widetilde{e}_{N}^{i}$. From the supporting hyperplane theorem \cite{Rockafellar}, $\varphi=1$ (or $\varphi=-1$) and $\widetilde{\varphi}=a^{j}$ hold for some $a^{j}\neq \pm \textbf{1}_{N}$, where $a^{j}$ is defined in (\ref{eq8}). It yields $\langle a^{j}, v\rangle\in(-\infty,1)\cup(1,\infty)$ and $\langle a^{j}, \widetilde{e}_{N}^{i}\rangle\in[-1,1]$. According to the separating hyperplane theorem \cite{Rockafellar}, one has $v\notin\mathcal{P}_\mathcal{V}$, which is contradiction.

\emph{Case II}: $\varphi$ is not the largest signed distance along the direction $\widetilde{\varphi} $ from the origin of $ \mathcal{P}_\mathcal{V}$ to a hyperplane $\mathcal{H}'$ containing $\widetilde{e}_{N}^{i}$. According to the separating hyperplane theorem \cite{Rockafellar}, one can design $\mathcal{H}'$ under $ \widetilde{\varphi}$ and $\varphi\in(-1,1)$, such that $\mathcal{H}'$ separates $ \mathcal{P}_\mathcal{V}$ into $\mathcal{P}_\mathcal{V}\backslash \{v\}$ and $\{v\}$. For any $x\in\mathcal{P}_\mathcal{V}\backslash \{v\}$, $x$ can be always written by a form of a convex combination of the entries in $\mathcal{V}$. This is a contradiction, since the set of all extreme points is the smallest subset of $\mathcal{P}_\mathcal{V}$ whose convex hull is equal to $\mathcal{P}_\mathcal{V}$ \cite{Lay}. Thus, $v$ is not an extreme point.

Based on the above discussions, $\text{ext}\mathcal{P}_\mathcal{V}=\mathcal{V}$. The proof is completed.
\end{IEEEproof}

\section{Proof of Lemma \ref{lemma2}}\label{appendix2}
\begin{IEEEproof}
The proof is equivalent to show $\overline{\mathcal{P}}_{2N}=\mathcal{P}_{\overline{\mathcal{V}}}$ and $\text{ext}\mathcal{P}_{\overline{\mathcal{V}}}=\overline{\mathcal{V}}$. For the case of $\text{ext}\mathcal{P}_{\overline{\mathcal{V}}}=\overline{\mathcal{V}}$, it is similar with Step II in the proof of Lemma \ref{lemma1}, which is omitted. The following proof only considers $\overline{\mathcal{P}}_{2N}=\mathcal{P}_{\overline{\mathcal{V}}}$.

For any $\overline{x}\in\overline{\mathcal{P}}_{2N}$, $\langle a^{i},\overline{x}\rangle\in[-1,1]$ and $\langle \textbf{1}_{N}\otimes e_{2}^{j},\overline{x}\rangle=0$ hold for all $i=1,\ldots,2^{N}$ and $j=1,2$, where $a^{i}\in\mathbb{R}^{2N}$ has the same definition of its counterpart in (\ref{eq8}). Let $\overline{x}=\overline{x}^{1}\otimes e_{2}^{1}+ \overline{x}^{2}\otimes e_{2}^{2}$, and $\overline{x}^{1},\overline{x}^{2}\in\mathbb{R}^{N}$, there must be $\|\overline{x}^{1} \|_{1}\leq1$ and $\|\overline{x}^{2} \|_{1}\leq1-\|\overline{x}^{1} \|_{1}$ from the equivalent form of $\mathcal{P}_{N}$ in (\ref{eq9}). From Lemma \ref{lemma1}, one has $\overline{x}^{1}\in\|\overline{x}^{1} \|_{1}\mathcal{P}_{ \mathcal{V} }$ and $\overline{x}^{1}\in(1-\|\overline{x}^{1} \|_{1})\mathcal{P}_{ \mathcal{V} }$, which means that there exist two nonnegative sequences $\{\overline{\lambda}_{j}^{i}\}_{i=1}^{\widetilde{N}}$ $(j=1,2)$ satisfying $\sum_{i=1}^{\widetilde{N}}\lambda_{1}^{i}=\|\overline{x}^{1} \|_{1}$, $\sum_{i=1}^{\widetilde{N}}\lambda_{2}^{i}=\|\overline{x}^{2} \|_{1}\leq1-\|\overline{x}^{1} \|_{1}$, and $\lambda_{j}^{i}\in[0,1]$, such that
\begin{align}
\overline{x}=&\sum_{p=1}^{\widetilde{N}} \lambda_{2}^{p}\Big(\frac{\overline{x}^{1}\otimes e^{1}_{2}}{\|\overline{x}^{2} \|_{1}}+\widetilde{e}_{N}^{p}\otimes e^{2}_{2}\Big)\nonumber \\
=& \sum_{p=1}^{\widetilde{N}} \lambda_{2}^{p}\sum_{q=1}^{\widetilde{N}} \lambda_{1}^{q}\Big(\frac{\widetilde{e}_{N}^{q}\otimes e^{1}_{2}}{\|\overline{x}^{2} \|_{1}}+\frac{\overline{e}_{N,2}^{p}}{\|\overline{x}^{1} \|_{1}}\Big)\nonumber \\
=&\sum_{p,q}^{\widetilde{N}} \frac{\lambda_{2}^{p} \lambda_{1}^{q}}{\|\overline{x}^{1} \|_{1}\|\overline{x}^{2} \|_{1}}(\|\overline{x}^{1} \|_{1}\overline{e}_{N,1}^{q}+\|\overline{x}^{2} \|_{1}\overline{e}_{N,2}^{p})\nonumber \\
\leq& \sum_{p,q}^{\widetilde{N}} \frac{\lambda_{2}^{p} \lambda_{1}^{q}}{\|\overline{x}^{1} \|_{1}\|\overline{x}^{2} \|_{1}}(\|\overline{x}^{1} \|_{1}\overline{e}_{N,1}^{q}+(1-\|\overline{x}^{1} \|_{1})\overline{e}_{N,2}^{p}),\label{eql21}
\end{align}
holds for $\|\overline{x}^{1} \|_{1}\neq0$ and $\|\overline{x}^{2} \|_{1}\neq0$. It is easy to check that $\sum_{p,q}^{\widetilde{N}} \lambda_{2}^{p} \lambda_{1}^{q}/(\|\overline{x}^{1} \|_{1}\|\overline{x}^{2} \|_{1})=1$, which further implies that $\overline{x}$ can be written by a convex combination of entries in $\overline{\mathcal{V}}$. For $\|\overline{x}^{1} \|_{1}=0$ or $\|\overline{x}^{2} \|_{1}=0$, obviously, $\overline{x}$ can be written by a form of a convex combination of entries in $\overline{\mathcal{V}}$. Therefore, $\overline{\mathcal{P}}_{2N}\subseteq\mathcal{P}_{\overline{\mathcal{V}}}$ holds.

For any $\overline{x}\in\mathcal{P}_{\overline{\mathcal{V}}}$, let $\overline{x}=\sum_{j=1}^{2}\sum_{i=1}^{\widetilde{N}}\lambda_{j}^{i}\overline{e}_{N,j}^{i}$, where $\lambda_{j}^{i}\in[0,1]$ and $\sum_{j=1}^{2}\sum_{i=1}^{\widetilde{N}}\lambda_{j}^{i}=1$. Then, for $p=1,\ldots,2^{N+1}$, one has
\begin{align}
\langle a^{p},\overline{x}\rangle=\sum_{j=1}^{2}\sum_{i=1}^{\widetilde{N}}\lambda_{j}^{i}\langle a^{p}, \overline{e}_{N,j}^{i}\rangle
\in[-1,1],\label{eql22}
\end{align}
where $a^{p}\in\mathbb{R}^{2N}$ has the same definition of its counterpart in (\ref{eq8}). It means that $\overline{x}\in\mathcal{P}_{2N}$ holds. Moreover, for $p=1,2$ and $j=1,2$, it is easy to get
\begin{align}
\langle\textbf{1}_{N}\otimes e_{2}^{p},\overline{x}\rangle=\sum_{j=1}^{2}\sum_{i=1}^{\widetilde{N}}\lambda_{j}^{i}\langle \textbf{1}_{N}\otimes e_{2}^{p}, \widetilde{e}_{N}^{i}\otimes e^{j}_{2}\rangle=0,\label{eql23}
\end{align}
which implies that $\overline{x}\in\overline{\mathcal{H}}$ holds. From (\ref{eql22}) and (\ref{eql23}), one has $\overline{x}\in\overline{\mathcal{P}}_{2N}$. The proof is completed.
\end{IEEEproof}

\section{Proof of Lemma \ref{lemma3}}\label{appendix3}
\begin{IEEEproof}
From Lemma \ref{lemma2}, one knows $\overline{\mathcal{P}}_{2N}=\mathcal{P}_{\overline{\mathcal{V}}}$. Then, for all $j=1,2$ and $k\in\mathbb{N}$, it only needs to show $\mathcal{S}_{k}^{T} \overline{e}_{N,j}^{i}\subset\mathcal{P}_{\overline{\mathcal{V}}}$.

For the case of $j=1$, one has
\begin{align}
\|\mathcal{S}_{k}^{T}  \overline{e}_{N,1}^{i}\|_{1}  \leq&\max_{p,q}\sum_{\imath=1,p\neq q}^{N}\frac{|s_{k}^{p\imath}-s_{k}^{q\imath}|}{2}\nonumber \\
=&1-\min_{p,q}\sum_{\imath=1}^{N}\min\{s_{k}^{p\imath},s_{k}^{q\imath}\}\label{eql31}.
\end{align}
The following proof is inspired by \cite{Liu2020}. For the matrix $\mathcal{S}_{k}\in\mathbb{S}_{\beta}$, let $S_{pq}^{1}=\{\imath|s^{p\imath}_{k}\leq s^{q\imath}_{k}\}$, $S_{pq}^{1+}=\{\imath|s^{p\imath}_{k}\geq0\}$, $S_{pq}^{1-}=S_{pq}^{1}\backslash S_{pq}^{1+}$, $S_{pq}^{2}=\{\imath|s^{p\imath}_{k}> s^{q\imath}_{k}\}$, $S_{pq}^{2+}=\{\imath|s^{q\imath}_{k}\geq0\}$, $S_{pq}^{2-}=S_{pq}^{2}\backslash S_{pq}^{2+}$, and $\imath_{mn}=\arg\max_{\imath}\{\min\{s^{p\imath}_{k},s^{q\imath}_{k}\}\}$ ($m\neq n$), then $\min\{s^{p\imath_{mn}}_{k},s^{q\imath_{mn}}_{k}\}\geq \alpha$. Since the graph is neighbor shared, it follows from (\ref{eql31}) that
\begin{align}
\sum_{\imath=1}^{N}\min\{s_{k}^{p\imath},s_{k}^{q\imath}\}=&\sum_{\imath\in S_{pq}^{1}}s_{k}^{p\imath}+\sum_{\imath\in S_{pq}^{2}}s_{k}^{q\imath}\nonumber \\
=&\min\{s^{p\imath_{mn}}_{k},s^{q\imath_{mn}}_{k}\}+\sum_{\imath\in S_{pq}^{1+}\backslash \{\imath_{mn}\}}s_{k}^{p\imath}\nonumber \\
&+\sum_{\imath\in S_{pq}^{1-}}s_{k}^{p\imath}+\sum_{i\in S_{pq}^{2+}\backslash \{\imath_{mn}\}}s_{k}^{q\imath}\nonumber \\
&+\sum_{\imath\in S_{pq}^{2-}}s_{k}^{q\imath}\nonumber \\
\geq&\alpha+\sum_{\imath\in S_{pq}^{1-}}s_{k}^{p\imath}+\sum_{i\in S_{pq}^{2-}}s_{k}^{q\imath}\nonumber \\
\geq&\alpha-2\beta.\label{eql32}
\end{align}
From (\ref{eql31}) and (\ref{eql32}), one gets
\begin{align}
\|\mathcal{S}_{k}^{T}  \overline{e}_{N,1}^{i}\|_{1}\leq&\eta(S_{k})\leq1-\alpha+2\beta<1.\label{eql33}
\end{align}
Moreover,  for all $i=1,\ldots,\widetilde{N}$ and $p=1,2$, one has $\langle\textbf{1}_{N}\otimes e_{2}^{p},\mathcal{S}_{k}^{T}\overline{e}_{N,1}^{i}\rangle=0$. Therefore, $\mathcal{S}_{k}^{T} \overline{e}_{N,1}^{i}\subset\mathcal{P}_{\overline{\mathcal{V}}}$ holds.

For the case of $j=2$, one can obtain the similar result from (\ref{eql31})-(\ref{eql33}). It completes the proof.
\end{IEEEproof}

\section{Proof of Lemma \ref{lemma4}}\label{appendix4}
\begin{IEEEproof}
For any $x\in\mathcal{P}_{\overline{\mathcal{V}}}$, one has $x=\sum_{j=1}^{2}\sum_{i=1}^{\widetilde{N}}\lambda_{j}^{i}\overline{e}_{N,j}^{i}$, where $\lambda_{j}^{i}\in[0,1]$ and $\sum_{j=1}^{2}\sum_{i=1}^{\widetilde{N}}\lambda_{j}^{i}=1$. From (\ref{eql31})-(\ref{eql33}), one has
\begin{align}
\|\mathcal{S}_{k}^{T}  x\|_{1}\leq&\eta(S_{k})\sum_{j=1}^{2}\sum_{i=1}^{\widetilde{N}}\lambda_{j}^{i}\leq 1-\alpha+2\beta<1.\label{eql41}
\end{align}
From Lemma \ref{lemma3}, one knows $x\in\mathcal{P}_{2N}$, which implies that $x$ has a form of convex combination of $\{\pm e_{2N}^{i}\}_{i=1}^{2N}$. Without loss of generality, let $x=\sum_{i=1}^{2N}(\lambda^{i}e_{2N}^{i}-\mu^{i}e_{2N}^{i})$, where $\lambda^{i},\mu^{i}\in[0,1]$, and $\sum_{i=1}^{2N}(\lambda^{i}+\mu^{i})=1$. Thus, one has
\begin{align}
\|\mathcal{S}_{k}^{T}  x\|_{1}=&\Big\|(S_{k}^{T}\otimes I_{2}) \sum_{i=1}^{2N}(\lambda^{i}e_{2N}^{i}-\mu^{i}e_{2N}^{i})\Big\|_{1}\nonumber\\
=&\Big\|\sum_{j=1}^{N}\sum_{i=1}^{ N} s^{ij}_{k}(\lambda^{2i-1}e_{2N}^{2i-1}-\mu^{2i-1}e_{2N}^{2i-1}+ \lambda^{2i}e_{2N}^{2i} \nonumber\\
&-\mu^{2i}e_{2N}^{2i})\Big\|_{1}\nonumber\\
\leq&\max_{i}\Big\{ \sum_{\substack{j=1\\s^{ji}_{k}\geq0}  }^{N}s^{ij}_{k}+\sum_{\substack{j=1\\s^{ji}_{k}<0}}^{N}|s^{ij}_{k}|\Big\}\sum_{i=1}^{2N}(\lambda^{i}+\mu^{i})\nonumber\\
\leq&1+2\beta.\label{eql42}
\end{align}
Then, $\mathcal{S}_{k}^{T}\mathcal{P}_{\overline{\mathcal{V}}}\subseteq \mathcal{S}_{k}^{T}\mathcal{P}_{2N}$ holds. Moreover, it is not hard to check $\mathcal{S}_{k}^{T}\mathcal{P}_{\overline{\mathcal{V}}}\subseteq \overline{\mathcal{H}}$. Therefore, $\mathcal{S}_{k}^{T}\mathcal{P}_{\overline{\mathcal{V}}}\subseteq\mathcal{S}_{k}^{T}\mathcal{P}_{2N}\cap\overline{\mathcal{H}}$.

For the converse, it can be regarded as an inclined projection $\mathcal{S}_{k}^{T}\mathcal{P}_{2N}$ onto $\overline{\mathcal{H}}$, which is similar with the proofs in Lemmas \ref{lemma1} and \ref{lemma3}. The details are omitted. Hence, the proof is completed.
\end{IEEEproof}

\section{Proof of Proposition \ref{proposition2}}\label{appendix5}

\begin{IEEEproof}
Obviously, $\mathcal{R}_{k}\mathcal{P}_{\overline{\mathcal{V}}}\subseteq \mathcal{R}_{k}\mathcal{P}_{2N}$ due to $ \mathcal{P}_{\overline{\mathcal{V}}}\subseteq  \mathcal{P}_{2N}$. The rest proof only considers $\mathcal{R}_{k}\mathcal{P}_{2N}\subseteq \mathcal{P}_{2N}^{\circ}$. From (\ref{eq8}), for any $x,y\in\mathcal{P}_{2N}$, $\|x\|_{1}\leq1$ and $\|y\|_{1}\leq1$, one has
\begin{align}
\langle\mathcal{R}_{k}x,y\rangle=&\sum_{i=1}^{N}\Big[(x^{2i-1}\cos\theta^{i}_{k}-x^{2i}\sin\theta^{i}_{k})y^{2i-1}\nonumber\\
&+(x^{2i-1}\sin\theta^{i}_{k}+x^{2i}\cos\theta^{i}_{k})y^{2i}\Big]\nonumber\\
\leq&\frac{1}{2}\sum_{i=1}^{N}\Big[(x^{2i-1}\cos\theta^{i}_{k}-x^{2i}\sin\theta^{i}_{k})^{2}+(y^{2i-1})^{2}\nonumber\\
&+(x^{2i-1}\sin\theta^{i}_{k}+x^{2i}\cos\theta^{i}_{k})^{2}+(y^{2i})^{2}\Big]\nonumber\\
\leq&\frac{1}{2}\sum_{i=1}^{2N}\Big[(x^{i})^{2}+(y^{i})^{2}\Big]\nonumber\\
\leq&\frac{1}{2}\Big[\Big(\sum_{i=1}^{2N}|x^{i}|\Big)^{2}+\Big(\sum_{i=1}^{2N}|y^{i}|\Big)^{2}\Big]\leq1,\label{eqp21}
\end{align}
where $x=[x^{1},\ldots,$ $x^{2N}]^{T}$ and $y=[y^{1},\ldots,y^{2N}]^{T}$. Therefore, $\mathcal{R}_{k}\mathcal{P}_{2N}\subseteq \mathcal{P}_{2N}^{\circ}$ holds. The proof is completed.
\end{IEEEproof}

\section{Proof of Proposition \ref{proposition3}}\label{appendix6}

\begin{IEEEproof}
Let $\mathcal{P}_{2N}^{*}=\{\chi\in\mathbb{R}^{2N}||\chi^{i}|\leq1/2N\}$, where $\chi^{i}$ is the $i$-th entry of $\chi$. One needs to prove $h_{k}\mathcal{P}_{2N}^{\circ}\subseteq\mathcal{P}_{2N}^{*}\subseteq \mathcal{P}_{2N}$. For any $x=[x^{1},\ldots,x^{2N}]^{T}\in\mathcal{P}_{2N}^{\circ}$ and $y
\in\{\pm e_{2N}^{i}\}\subset\text{ext} \mathcal{P}_{2N}\subseteq\mathcal{P}_{2N}$, one knows $\langle  x,y \rangle \leq1$, which implies that $h_{k} |x^{i}|\leq1/2N$ holds for all $i=1,\ldots,2N$. One further derives $h_{k} \mathcal{P}_{2N}^{\circ}\subseteq\mathcal{P}_{2N}^{*}$. Moreover, for any $\chi\in\mathcal{P}_{2N}^{*}$, one has $\sum_{i=1}^{2N}|\chi^{i}|\leq1$. Thus, $\mathcal{P}_{2N}^{*}\subseteq\mathcal{P}_{2N}$ holds, and therefore, $h_{k}\mathcal{P}_{2N}^{\circ}\subseteq\mathcal{P}_{2N}$. It completes the proof.
\end{IEEEproof}

\section{Proof of Proposition \ref{proposition4}}\label{appendix7}

\begin{IEEEproof}
For any $y\in\mathcal{P}_{2N}$ and $x\in\mathcal{P}_{\overline{\mathcal{V}}}$, $\|y\|_{1}\leq1$, $\|x\|_{1}\leq1$, and $\langle\textbf{1}_{N}\otimes e_{2}^{i},x\rangle=0$ hold, where $x=[x^{1},\ldots,$ $x^{2N}]^{T}$, $y=[y^{1},\ldots,y^{2N}]^{T}$, and $i=1,2$. Similar with (\ref{eqp21}), one has
\begin{align}
\langle\Lambda\mathcal{R}_{k}&x,y\rangle\nonumber\\
\leq&\frac{1}{2}\Big[\sum_{i=1}^{N}\Big((x^{2i-1}\cos\theta^{i}_{k}-x^{2i}\sin\theta^{i}_{k})^{2}+(x^{2i-1}\sin\theta^{i}_{k}\nonumber\\
&+x^{2i}\cos\theta^{i}_{k})^{2}\Big)+\frac{1}{N^2}\sum_{i=1,i\neq j}^{N}\sum_{j=1}^{N} \Big((y^{2j-1}-y^{2i-1})^{2}\nonumber\\
&+(y^{2j}-y^{2i})^{2}\Big)\Big]\nonumber\\
\leq&\frac{1}{2}\Big[\Big(\sum_{i=1}^{2N}|x^{i}|\Big)^{2}+\Big(\sum_{i=1}^{2N}|y^{i}|\Big)^{2}\Big]\leq1,\label{eqp41}
\end{align}
which means $ \Lambda\mathcal{R}_{k}x \in \mathcal{P}_{2N}^{\circ}$. Then, for $i=1,2$, one has $\langle\textbf{1}_{N}\otimes e_{2}^{i},\Lambda\mathcal{R}_{k}x\rangle=\langle\Lambda\textbf{1}_{N}\otimes e_{2}^{i},\mathcal{R}_{k}x\rangle =\langle\textbf{0}_{2N},\mathcal{R}_{k}x\rangle=0$. Thus, one gets $\Lambda\mathcal{R}_{k}\mathcal{P}_{\overline{\mathcal{V}}}\subseteq\mathcal{P}_{2N}^{\circ}\cap\overline{\mathcal{H}}$, which completes the proof.
\end{IEEEproof}
\end{appendices}

\end{document}